\documentclass[11pt,twoside]{article}

\usepackage{amscd}
\usepackage{epsfig}
\usepackage{latexsym}
\usepackage{amsfonts}
\usepackage{amsthm}

\def\ps@myheadings{\let\@mkboth\@gobbletw
\def\@oddhead{\hbox{}
\rightmark\hfil\eightrm\thepage}   
\def\@oddfoot{}\def\@evenhead{\eightrm\thepage\hfil
\leftmark\hbox{}}\def\@evenfoot{}
\def\sectionmark##1{}\def\subsectionmark##1{}}
\def\runninghead#1#2{\pagestyle{myheadings}
\markboth{\centerline{\footnotesize\rm{#1}}}
{\centerline{\footnotesize\rm{#2}}}}

\def\abstract#1{{
	\centering{\begin{minipage}{4.2in}\footnotesize\baselineskip=10pt
	{\sc Abstract.\/} #1\par 
	\end{minipage}}\par}}

\renewenvironment{thebibliography}[1]
	{\frenchspacing
	\rm\baselineskip=11pt
	 \begin{list}{[\arabic{enumi}]}
	{\usecounter{enumi}\setlength{\parsep}{0pt}
	 \setlength{\leftmargin 13.7pt}{\rightmargin 0pt} 
	 \setlength{\itemsep}{2pt} \settowidth
	{\labelwidth}{[#1]}\sloppy}}{\end{list}}
 
\DeclareSymbolFont{EulerScript}{U}{eus}{m}{n}
\DeclareSymbolFontAlphabet\mathscr{EulerScript}

\newtheoremstyle{mythm}{11pt}{11pt}{\it}{}{\sc}{.}{ }{}
\theoremstyle{mythm}
\newtheorem{prop}{Proposition}[section]
\newtheorem{lemma}[prop]{Lemma}

\newtheorem{thm}{Theorem}[section]

\newtheoremstyle{myex}{11pt}{11pt}{\rm}{}{\sc}{.}{ }{}
\theoremstyle{myex}
\newtheorem*{ex}{Example}

\newcounter{itemlistc}
\newenvironment{itemlist}
    	{\setcounter{itemlistc}{0}
	 \begin{list}{-}
	{\usecounter{itemlistc}
	 \setlength{\parsep}{0pt}
	 \setlength{\itemsep}{5pt}}}{\end{list}}

\newcommand{\Z}{{\mathbb Z}}
\newcommand{\Q}{{\mathbb Q}}
\newcommand{\R}{{\mathbb R}}
\newcommand{\C}{{\mathbb C}}
\newcommand{\lk}{{\ell k}}
\newcommand{\N}{{\mathscr N}}
\newcommand{\V}{{\mathscr V}}
\newcommand{\A}{{\mathscr A}}
\newcommand{\PP}{{\mathscr P}}
\newcommand{\E}{{\mathscr E}}
\newcommand{\LL}{{\mathscr L}}
\newcommand{\m}{(\underline{m})}
\newcommand{\fb}{\overline{F}}

\newcommand{\rk}{\hbox{rk}\,}
\newcommand{\Ho}{\hbox{H}}
\newcommand{\im}{\hbox{Im}\,}
\newcommand{\Ker}{\hbox{Ker}\,}
\newcommand{\Hom}{\hbox{Hom}}

\begin{document}

\renewcommand{\thefootnote}{\fnsymbol{footnote}}

\renewcommand{\section}[1] {\vspace{15pt}\addtocounter{section}{1} 
\centerline{\bf\thesection. #1}\par\vspace{5pt}}
\newcommand{\nonumsection}[1] {\vspace{12pt}\centerline{\bf #1}
	\par\vspace{5pt}}
	
\runninghead{THE ALEXANDER MODULE OF LINKS AT INFINITY}{DAVID CIMASONI}

\thispagestyle{empty}

\centerline{\bf THE ALEXANDER MODULE OF LINKS AT INFINITY}
\vskip1cm
\centerline{\sc david cimasoni 
}
\vskip1cm
\abstract{Walter Neumann \cite{Neu} showed that the topology of a ``regular'' algebraic curve $V\subset\C^2$ is determined up to
proper isotopy by some link in $S^3$ called the link at infinity of $V$. In this note, we compute the Alexander module over
$\C[t^{\pm 1}]$ of any such link at infinity.}
\vskip0.5cm

\section{Introduction}
The intersection of a reduced algebraic curve $V\subset\C^2$ with any sufficiently large sphere $S^3$ about the origin in $\C^2$
gives a well-defined link called the {\em link at infinity\/} of $V\subset\C^2$. This link at infinity was first introduced by Walter
Neumann and Lee Rudolph \cite{N-R} and studied further by Neumann \cite{Neu}. In order to state several of their results, let us recall
some terminology. The fiber $f^{-1}(c)$ of a polynomial map $f\colon\C^2\to\C$ is called {\em regular\/} if there exists a neighborhood
$D$ of $c$ in $\C$ such that $f|\colon f^{-1}(D)\to D$ is a locally trivial fibration. An algebraic curve $V\subset\C^2$ is {\em regular}
if it is a regular fiber of its defining polynomial. One might think that if $c$ is
not a singular value of $f$, then $f^{-1}(c)$ is regular; this is wrong.
In fact, the following additional condition is required: a fiber $f^{-1}(c)$ is
{\em regular at infinity\/} if there exists a neighborhood $D$ of $c$ in $\C$ and a compact $K$ in $\C^2$ such that $f$ restricted to
$f^{-1}(D)\setminus K$ is a locally trivial fibration. It can be proved that $f^{-1}(c)$ is regular if and only if it is non-singular
and regular at infinity \cite{H-L}.

A first interesting result is that only finitely many fibers of a given $f$ are irregular at infinity, and that the regular fibers
all define the same link at infinity up to isotopy: it is called the {\em regular link at infinity\/} of $f$, and denoted by
$\LL(f,\infty)$. Furthermore, $\LL(f,\infty)$ is a fibered link if and only if all the fibers of $f$ are regular at infinity.
Finally, Walter Neumann proved the following striking result: the topology of a regular algebraic
curve $V\subset\C^2$ (as an embedded smooth manifold) is determined by its link at infinity. More precisely: up to isotopy in $S^3$,
there exists a unique minimal Seifert surface $F$ for $\LL(f,\infty)$, and $V$ is properly isotopic to the embedded surface obtained
from $F$ by attaching a collar out to infinity in $\C^2$ to the boundary of $F$.

In the present note, we give a closed formula for the Alexander module over $\C[t^{\pm 1}]$ of the regular link at infinity of any
polynomial map $f\colon\C^2\to\C$ (Theorem \ref{thm:infinity}). The decisive property of $\LL(f,\infty)$ is that it can be seen as
the boundary of the fiber $F$ of a fibered multilink (see \cite[Theorem 4]{Neu}). Furthermore, this multilink can be 
constructed by iterated cabling and connected sum operations from the unknot, and the
Alexander module over $\C[t^{\pm 1}]$ of this type of fibered multilinks is well-known (we recall this result of \cite{E-N} in Theorem
\ref{thm:C-fibered} below). Therefore, our method will be to consider a fibered multilink with fiber $F$ and Alexander module $A$,
and to compute the Alexander module of the oriented link $\LL=\partial F$ from the module $A$ (Proposition \ref{prop:L'}). This is
achieved by introducing ``generalized Seifert forms'' for the multilink, and comparing them with the traditional Seifert form for $\LL$.
The result is then applied to $\LL(f,\infty)$ (Theorem \ref{thm:infinity}), and an example concludes the paper. 

\section{Fibered multilinks}
A {\em multilink\/} \cite{E-N} is an oriented link $L=L_1\cup\dots\cup L_n$ in $S^3$
together with an integer $m_i$ associated with each component $L_i$, with the convention that a component $L_i$ with multiplicity
$m_i$ is the same as $-L_i$ ($L_i$ with reversed orientation) with mutliplicity $-m_i$. Throughout this paper, we will write
$\underline{m}$ for the integers $(m_1,\dots,m_n)$, $d$ for their greatest common divisor, and $L\m$ for the multilink.
Of course, a set of multiplicities $\underline{m}$ can be thought of as an element of $\Ho_1(L)$. If $X$ denotes the exterior of $L$, 
several classical theorems implie that $\Ho_1(L)$ is isomorphic to $[X,S^1]$, the group of homotopy classes of maps $X\to S^1$.
As a consequence, assigning a set of multiplicities to an oriented link is a way to specify a preferred infinite cyclic covering 
$\widetilde{X}\m\stackrel{p}{\to}X$: it is the pullback $\Z$-bundle $\phi^{\ast}\exp$, where $\R\stackrel{\exp}{\to}S^1$ is the 
universal $\Z$-bundle and $X\stackrel{\phi}{\to}S^1$ any map in the homotopy class $\underline{m}$: 
$$
\begin{CD}
\widetilde{X}\m @>>>\R\\
@V{p}VV	 @VV{\exp}V\\
X@>{\phi}>>\phantom{.}S^1.
\end{CD}
$$
Choosing a generator $t$ of the infinite cyclic group of the covering endows
$\Ho_\ast(\widetilde{X}\m)$ with a structure of module over $\Z\!<t>\,\,=\Z[t^{\pm 1}]$, the ring of Laurent polynomials with integer coefficients.
Most of these invariants are not interesting:
it is easy to prove that $\Ho_0(\widetilde{X}\m)\simeq\Z[t^{\pm 1}]/(t^d-1)$, that $\Ho_2(\widetilde{X}\m)$ is a free module with
the same rank as $\Ho_1(\widetilde{X}\m)$, and of course, that
$\Ho_i(\widetilde{X}\m)=0$ for all $i\ge 3$. Therefore, the only interesting module is $\Ho_1(\widetilde{X}\m)$: it is called the
{\em Alexander module\/} of the multilink $L\m$, and we will denote it by $A(L\m)$. Also, we will write $A(L\m;{\mathbb K})$ for the ${\mathbb K}[t^{\pm 1}]$-module
$A(L\m)\otimes{\mathbb K}[t^{\pm 1}]$, where ${\mathbb K}=\Q$ or $\C$.
Given $\PP$ an $m\times n$ presentation matrix of $A(L\m)$ (that is, the matrix corresponding to a finite presentation of $A(L\m)$ with $n$ generators
and $m$ relations), the greatest common divisor of the $n\times n$ minor determinants of ${\mathscr P}$ is called the {\em Alexander polynomial\/} of $L\m$.
This Laurent polynomial, denoted by $\Delta_{L\m}$, is only defined up to multiplication by units of $\Z[t^{\pm 1}]$.
Of course, if a multilink has multiplicities $\pm 1$, it is just an oriented link and these Alexander invariants
coincide with the usual Alexander invariants of the corresponding oriented link.

Let us now recall the definition of a very interesting class of multilinks that generalizes the notion of fibered link: a
{\em fibered multilink\/} is a multilink $L\m$ such that there exists a locally
trivial fibration $X\stackrel{\varphi}{\to}S^1$ in the homotopy class $\underline{m}\in[X,S^1]$. The oriented surface $F=\varphi^{-1}(1)$ is called the
{\em fiber\/} of $L\m$. The diagram
$$
\begin{CD}
\widetilde{X}\m @>{\Phi}>>\R\\
@V{p}VV	 @VV{\exp}V\\
X@>{\varphi}>>S^1
\end{CD}
$$
can now be understood as defining the pullback fibration $\Phi=\exp^\ast \varphi$. Since $\R$ is contractible, there exists a
homeomorphism $F\times\R\to\widetilde{X}\m$ such that the following diagram commutes:
$$
\begin{CD}
F\times\R@>>>\widetilde{X}\m\\ 
@V{\pi}VV  @VV{\Phi}V \\
\R @=\phantom{.}\R.	
\end{CD}
$$
Hence, the generator $\widetilde{X}\m\stackrel{t}{\to}\widetilde{X}\m$ of the infinite cyclic group of the covering $p$ can be
seen as the transformation
$$
\begin{array}{ccl}
F\times\R&\longrightarrow&F\times\R\cr
(x,z)&\longmapsto&(h(x),z+1),
\end{array}
$$
where $F\stackrel{h}{\longrightarrow}F$ is some homeomorphism, unique up to isotopy, called the {\em monodromy\/} of the multilink $L\m$.
We will use the same terminology for the induced automorphism $\Ho_1(F)\stackrel{h_\ast}{\longrightarrow}\Ho_1(F)$.

\begin{prop}\label{prop:H-tI}
A presentation matrix of the Alexander module of a fibered multilink is given by $H^T-tI$, where $H$ is any matrix of the monodromy.
In particular, the Alexander polynomial of a fibered multilink is the characteristic polynomial of the monodromy.
\end{prop}
\proof As seen in the above discussion, there is an isomorphism of $\Z$-modules $\Ho_1(F)\stackrel{f}{\to}\Ho_1(\widetilde{X}\m)$ such that
$t\cdot f(x)=f(h_\ast(x))$. Choosing a $\Z$-basis $e_1,\dots,e_\mu$ of $\Ho_1(F)$, this gives an exact sequence of $\Z[t^{\pm 1}]$-modules
$$
\bigoplus_{i=1}^\mu \Z[t^{\pm 1}]\,e_i\stackrel{h_\ast-t}{\longrightarrow}\bigoplus_{i=1}^\mu \Z[t^{\pm 1}]\,e_i\stackrel{f_\ast}
{\longrightarrow}\Ho_1(\widetilde{X}\m)\longrightarrow 0,
$$
where $f_\ast$ denotes the $\Z[t^{\pm 1}]$-linear extension of $f$. This is a finite presentation of $\Ho_1(\widetilde{X}\m)$, so 
$(H-tI)^T$ is a presentation matrix of this module. \endproof

\begin{figure}[Htb]
   \begin{center}
     \epsfig{figure=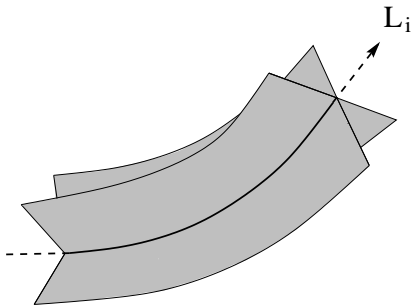,height=4cm}
     \caption{\footnotesize $\fb$ near the component $L_i$ of a multilink, with $m_i=4$.}
     \label{fig:fig1}
   \end{center}
\end{figure}

Let $F\subset S^3\setminus L$ be the fiber of a fibered multilink $L\m$, and let us denote by $\fb$ the union $F\cup L$ (see Figure \ref{fig:fig1}
for an illustration of $F$ and $\fb$ near a component of the multilink). The {\em Seifert forms\/} associated to $F$ are the bilinear forms
$$
\alpha_+,\alpha_-\colon\Ho_1(F)\times\Ho_1(\fb)\longrightarrow\Z
$$
given by $\alpha_+(x,y)=\lk(i_+x,y)$ and $\alpha_-(x,y)=\lk(i_-x,y)$, where $\lk$ denotes the linking number and
$i_+,i_-\colon\Ho_1(F)\to\Ho_1(S^3\setminus\fb)$ the morphisms induced by the push in the positive or negative normal direction off $F$.
We will use the notation $V_+,V_-$ for matrices of these forms.

As in the usual case of a fibered oriented link, the monodromy can be recovered from
the Seifert forms.

\begin{prop}\label{prop:fiber}
If a multilink is fibered with fiber $F$, the matrices $V_+$ and $V_-$ are square and unimodular. 
Furthermore, a matrix of the monodromy is given by $H=(V_+V_-^{-1})^T$.
\end{prop}
\proof The bilinear form $\alpha_+\colon\Ho_1(F)\times\Ho_1(\fb)\to\Z$ can be understood as a homomorphism
$\Ho_1(F)\to\Hom(\Ho_1(\fb),\Z)\simeq\Ho^1(\fb)$.
The composition of this morphism with the Alexander isomorphism $\Ho^1(\fb)\simeq\Ho_1(S^3\setminus\fb)$ is nothing but
$i_+\colon\Ho_1(F)\to\Ho_1(S^3\setminus\fb)$. The same holds for $\alpha_-$ and $i_-$. As a consequence, the Seifert matrix $V_+$ (resp. $V_-$) with respect
to basis $\A$ of $\Ho_1(F)$ and $\overline\A$ of $\Ho_1(\fb)$ is equal to the transposed matrix of
$i_+$ (resp. $i_-$) with respect to the basis $\A$ and $\overline\A^\ast$,
where $\overline\A^\ast$ is the dual basis of $\overline\A$ via Alexander duality.

Now, the fibration $S^3\setminus L\to S^1$
yields a fibration $S^3\setminus\fb\to (0,1)$, so $S^3\setminus\fb$ is homeomorphic to $F\times (0,1)$. Hence, the maps
$i_+,i_-\colon F\to S^3\setminus\fb$ are homotopy equivalences, and  $i_+,i_-\colon\Ho_1(F)\to\Ho_1(S^3\setminus\fb)$ are isomorphisms.
Therefore, the matrices $V_+,V_-$ are unimodular. Finally, the monodromy of a fibered multilink can be defined as the composition
$(i_-)^{-1} \circ (i_+)$. Therefore, a matrix of the monodromy is given by $H=(V_-^{-1})^TV_+^T=(V_+V_-^{-1})^T$.
\endproof
As an immediate consequence of this proposition, $\Ho_1(F)$ and $\Ho_1(\fb)$ have the same rank. We need some more information about these modules.

\begin{lemma}\label{lemma:homology}
Let $L\m$ be a fibered multilink with fiber $F$ of genus $g$. For $i=1,\dots,n$, $F$ has
$d_i=\gcd(m_i,\sum_{j\neq i}m_j\lk(L_i,L_j))$ boundary components near $L_i$. Furthermore, the homology of $F$ has the form
$$
{\mathrm H}_1(F)\;=\;G\;\oplus\bigoplus_{i=1\dots n-1\atop j=1\dots d_i}\Z T_i^j\;\oplus\;\bigoplus_{j=1}^{d_n-d}\Z T_n^j\,,
$$
where $G$ is a free $\Z$-module of rank $2dg$, and $T_i^1,\dots,T_i^{d_i}$ are the boundary components of $F$ near $L_i$.
Finally,
$$
{\mathrm H}_1(\fb)\;=\;G\;\oplus\left(\bigoplus_{i=1}^n\Z L_i\Bigg/\sum_{i=1}^n\frac{m_i}{d}L_i\right)\oplus\;\overline{B}\,,
$$
where $\overline{B}$ is a free $\Z$-module of rank $1-n-d+\sum_{i=1}^nd_i$.
\end{lemma}
\proof The fact that $F\cap\N(L_i)$ is a link with $d_i$ components is very easy to check and well-known (see \cite[p. 30]{E-N}).
Since $F$ consists of $d$ parallel copies of the fiber of the multilink $L(\frac{\underline{m}}{d})$, it may be assumed that $d=1$.
In this case, $F$ is a connected oriented surface of genus $g$ with $\sum_{i=1}^nd_i$ boundary components and the result holds.

We will now compute $\Ho_1(\fb)$ by induction on $d\ge 1$. Let us assume that $d=1$. The Mayer-Vietoris exact sequence associated
with the decomposition $\fb=F\cup(\fb\cap\N(L))$ gives
$$
0\to\Ho_1(\partial F)\stackrel{\varphi_1}{\to}\Ho_1(F)\oplus\Ho_1(L)\to\Ho_1(\fb)\to\widetilde{\Ho}_0(\partial F)
\stackrel{\varphi_0}{\to}\widetilde{\Ho}_0(L),
$$
where $\varphi_1(T_i^j)=(T_i^j,\frac{m_i}{d_i}L_i)$. Using the value of $\Ho_1(F)$, it follows that
$(\Ho_1(F)\oplus\Ho_1(L))/\im\varphi_1=G\oplus\left(\bigoplus_{i=1}^n\Z L_i\big/\sum_i m_i L_i\right)$.
Since the module $\Ker\varphi_0$ is free of rank $\sum_{i=1}^n(d_i-1)$, this concludes the case $d=1$.
Let us now consider a fibered multilink $L\m$ with $\gcd(m_1,\dots,m_n)=d>1$. Clearly, $\fb=\fb'\cup\fb''$, where
$F'$ (resp. $F''$) is the fiber of $L(\frac{\underline{m}}{d})$ (resp. $L(\frac{d-1}{d}\,\underline{m})$). The associated Mayer-Vietoris
sequence together with the case $d=1$ and the induction hypothesis give the result. \endproof

\begin{prop}\label{prop:split}
Let $L\m$ be a fibered multilink. For $i=1,\dots,n$, let us note $D_i=\gcd(d_1,\dots,d_i)$ with $d_i$ as above.
Then, the Alexander module of $L\m$ naturally factors into $A(L\m)\;=\;A_G\oplus A_B$, where
$$
A_B=\bigoplus_{i=1}^{n-1}\Z[t^{\pm 1}]\Big/\Bigl(\frac{(t^{D_i}-1)(t^{d_{i+1}}-1)}{(t^{D_{i+1}}-1)}\Bigr).
$$
\end{prop}
\proof As seen above, the fiber $F$ is given by $d$ parallel copies of a connected surface $\widetilde{F}$ with $\sum_{i=1}^n\frac{d_i}{d}$
boundary components. Let us write $\widetilde{F}=\widetilde{G}\cup\widetilde{B}$, where $\widetilde{G}$ is a closed surface with a single boundary component,
and $\widetilde{B}$ a planar surface with $1+\sum_{i=1}^n\frac{d_i}{d}$ boundary components. The Mayer-Vietoris sequence gives
$\Ho_1(\widetilde{F})=\Ho_1(\widetilde{G})\oplus\Ho_1(\widetilde{B})$. Therefore, $\Ho_1(F)=\Ho_1(G)\oplus\Ho_1(B)$, where $G$ (resp. $B$) consists of
$d$ parallel copies of $\widetilde{G}$ (resp. $\widetilde{B}$).
Since the monodromy $F\stackrel{h}{\to}F$ of $L\m$ is a homeomorphism, the monodromy $\Ho_1(F)\stackrel{h_\ast}{\to}\Ho_1(F)$
splits into $h_G\oplus h_B$, where $h_G=(h|_G)_\ast$ and $h_B=(h|_B)_\ast$. Therefore, a matrix $H$ of $h_\ast$ with respect to some
basis $\A=\A_G\cup\A_B$ of $\Ho_1(F)=\Ho_1(G)\oplus \Ho_1(B)$ can be written $H=H_G\oplus H_B$. By Proposition \ref{prop:H-tI}, $A(L\m)$ is presented by
$$
H^T-tI\;=\;H_G^T\oplus H_B^T-tI\;=\;(H_G^T-tI)\oplus(H_B^T-tI).
$$
Let us denote by $A_G$ (resp. $A_B$) the $\Z[t^{\pm 1}]$-module presented by $H_G^T-tI$ (resp. $H_B^T-tI$).
It remains to compute the module $A_B$.

As seen in Lemma \ref{lemma:homology}, a basis of $\Ho_1(B)$ is given by
$$
\A_B=\left<T_1^1,\dots,T_1^{d_1},\dots,T_{n-1}^1,\dots,T_{n-1}^{d_{n-1}},T_n^1,\dots,T_n^{d_n-d}\right>,
$$
where $T_i^1,\dots,T_i^{d_i}$ are the boundary components of $F$ near $L_i$. Clearly, the monodromy cyclically permutes these components, that is:
$h_\ast(T_i^j)=T_i^{j+1}$ for $1\le j\le d_i-1$ and $h_\ast(T_i^{d_i})=T_i^1$. Note that 
$$
h_\ast(T_n^{d_n-d})=T_n^{d_n-d+1}=-\sum_{i=1\dots n-1\atop j\equiv 1\!\!\!\!\!\pmod d}T_i^j - \sum_{1\le j\le d_n-d\atop j\equiv 1\!\!\!\!\!\pmod d}T_n^j
$$
in $\Ho_1(B)$, since $\partial\widetilde{F}=\bigcup_{i=1\dots n}\bigcup_{j\equiv 1\!\!\!\!\!\pmod d}T_i^j$. Therefore,
the matrix of $h_B$ with respect to $\A_B$ is given by
$$
H_B=\pmatrix{P_1& & & &\cr
		&\ddots& & & v\cr
		&	&P_{n-1}& &\cr
		&	&	&Q_n&},
$$ 
where $P_i$ is the $d_i\times d_i$-matrix
$$
P_i=\pmatrix{&&&1\cr 1&&&\cr &\ddots&&\cr &&1&},
$$
$Q_n$ the $(d_n-d)\times(d_n-d-1)$-matrix
$$
Q_n=\pmatrix{0&\dots&0\cr 1&&\cr &\ddots&\cr &&1\cr},
$$ 
and $v=(v_j)$ a vector such that $v_j=-1$ if $j\equiv 1\!\!\!\pmod d$, $v_j=0$ else.
It is easy to show that $H_B^T-tI$ is equivalent to
$$
\pmatrix{t^{d_1}-1\cr
	&	t^{d_2}-1\cr
	&	&	\ddots\cr
	&	&	&	t^{d_{n-1}}-1\cr
	\frac{t^{d_1}-1}{t^d-1}&\frac{t^{d_2}-1}{t^d-1}&\dots&\frac{t^{d_{n-1}}-1}{t^d-1}&\frac{t^{d_n}-1}{t^d-1}}
$$
as a presentation matrix. It is then an exercise to check that the module $A_B$ presented by $H^T_B-tI$ is equal to
$\bigoplus_{i=1}^{n-1}\Z[t^{\pm 1}]\Big/\Bigl(\frac{(t^{D_i}-1)(t^{d_{i+1}}-1)}{(t^{D_{i+1}}-1)}\Bigr)$.\endproof

We are finally ready to prove the main result of this paragraph.

\begin{prop}\label{prop:L'}
Let $L\m$ be a fibered multilink with fiber $F$ and multiplicities $m_i\neq 0$ for all $i$, and let us denote by $L'$ the oriented link given by
the boundary $\partial F$ of $F$. Then, the Alexander module of $L'$ over $\Q[t^{\pm 1}]$ is given by
$$
A(L';\Q)\;=\;(A_G\otimes\Q[t^{\pm 1}])\oplus\left(\Q[t^{\pm 1}]/(t-1)\right)^{n-1}\oplus\left(\Q[t^{\pm 1}]\right)^{\sum_{i=1}^n(d_i-1)},
$$
where $A_G$ is the direct factor of $A(L\m)$ given in Proposition \ref{prop:split}, and $d_i=\gcd(m_i,\sum_{j\neq i}m_j\lk(L_i,L_j))$.
\end{prop}
\proof Since $m_i\neq 0$, it follows that $d_i\neq 0$ for all $i$. By Lemma \ref{lemma:homology}, one can write
\begin{eqnarray*}
\Ho_1(F;\Q)&=&(G\otimes\Q)\;\oplus\;\bigoplus_{i=1}^{n-1}\Q\,\sum_{j=1}^{d_i}T_i^j
\;\oplus\bigoplus_{i=1\dots n-1\atop j=1\dots d_i-1}\Q T_i^j\;\oplus\;\bigoplus_{j=1}^{d_n-d}\Q T_n^j\\
\Ho_1(\fb;\Q)&=&(G\otimes\Q)\;\oplus\;\bigoplus_{i=1}^{n-1}\Q\,m_iL_i\;\oplus\;(\overline{B}\otimes\Q).
\end{eqnarray*}
The matrices $V_+$ and $V_-$ with respect to these basis of $\Ho_1(F;\Q)$ and $\Ho_1(\fb;\Q)$ are of the form
$$
V_+\;=\bordermatrix{
&\overbrace{\quad}^{2dg}&\overbrace{\quad}^{n-1}&\phantom{\overbrace{\quad}}\cr
\hfill\scriptstyle 2dg\Bigl\{&N&M^T&\ast\cr
\hfill\scriptstyle n-1\Bigl\{&M&\ell&\ast\cr
\hfill\phantom{\Bigl\{}&\ast&\ast&\ast\cr}
\qquad
V_-\;=\bordermatrix{
&\overbrace{\quad}^{2dg}&\overbrace{\quad}^{n-1}&\phantom{\overbrace{\quad}}\cr
\hfill\scriptstyle 2dg\Bigl\{&N^T&M^T&\ast\cr
\hfill\scriptstyle n-1\Bigl\{&M&\ell^T&\ast\cr
\hfill\phantom{\Bigl\{}&\ast&\ast&\ast\cr}.
$$
As seen in the proof of Proposition \ref{prop:split}, a matrix $H$ of the monodromy splits into $H_G\oplus H_B$. Furthermore, the basis of $\Ho_1(F;\Q)$
was chosen such that
$$
H_B=\pmatrix{I_{n-1}&\ast\cr 0&\ast\cr},
$$
where $I_{n-1}$ denotes the identity matrix of dimension $n-1$. By Proposition \ref{prop:fiber}, $V_+=H^TV_-$, that is,
$$
\scriptstyle\pmatrix{N&M^T&\ast\cr M&\ell&\ast\cr \ast&\ast&\ast\cr}=
\pmatrix{H_G^T&0&0\cr 0&I_{n-1}&0\cr 0&\ast&\ast\cr}\pmatrix{N^T&M^T&\ast\cr M&\ell^T&\ast\cr\ast&\ast&\ast\cr}=
\pmatrix{H_G^TN^T&H_G^TM^T&\ast\cr M&\ell^T&\ast\cr \ast&\ast&\ast\cr}.
$$
Therefore, we have the equalities
$$
N^T=NH_G\;,\quad M=MH_G\,,\quad\hbox{and}\quad\ell=\ell^T.\eqno{(\ast)}
$$

Let us keep them in mind, and turn to the computation of the Alexander module of $L'$. Since $F$ has $d$ connected components,
a connected Seifert surface $F'$ for $L'$ is obtained from $F$ via $d-1$ handle attachments.
Since $d_i\neq 0$ for all $i$, we can write
$$
\Ho_1(F';\Q)\;\,=\,\;(G\otimes\Q)\;\oplus\bigoplus_{i=1\dots n-1\atop j=1\dots d_i}\Q(d_iT_i^j)\;\oplus\;\bigoplus_{j=1}^{d_n-1}\Q(d_nT_n^j)\,.
$$
The Seifert matrix of $L'$ with respect to this basis has the form
$$
V'=\bordermatrix{&&&\cr
\hfill\scriptstyle 2dg\bigl\{&N&\ast&\ast\cr
\hfill\scriptstyle d_1\bigl\{&&\widetilde\ell_1&\cr
&&\vdots&\cr
\hfill\scriptstyle d_n-1\bigl\{&&\widetilde\ell_n&\cr}=
\bordermatrix{&\overbrace{\quad}^{2dg}&\overbrace{\quad}^{d_1}&&\overbrace{\quad}^{d_n-1}\cr
&N&&&&\cr
&\ast&{\widetilde{\ell_1}}^T&\dots&{\widetilde{\ell_n}}^T\cr
&\ast&&&&\cr},
$$ 
where $\widetilde\ell_i$ denotes $d_i$ copies of the same line $\ell_i$ ($d_n-1$ copies if $i=n$). A presentation matrix of
$A(L';\Q)$ is given by $\PP'=V'-t(V')^T$. Since $\sum_{i,j}T^j_i=0$ in $\Ho_1(F')$, it follows that $\ell_n=-\sum_{i=1}^{n-1}\ell_i$.
As a presentation matrix, $\PP'$ is therefore equivalent to
$$
\scriptstyle\pmatrix{N-tN^T&\ast&\ast\cr
&\ell_1(1-t)&\cr
&\vdots&\cr
&\ell_{n-1}(1-t)&\cr}=
\pmatrix{N-tN^T&&\cr\ast&\ell_1^T(1-t)\dots\ell_{n-1}^T(1-t)&0\dots 0\cr\ast&&\cr},
$$
where the number of zero columns is equal to 
$$
\sum_{i=1}^{n-1}(d_i-1)+(d_n-1)=\sum_{i=1}^n(d_i-1).
$$
With the notations used above for $V_+$ and $V_-$, this matrix is nothing but
$$
\pmatrix{N-tN^T&M^T(1-t)&0\dots 0\cr M(1-t)&\ell(1-t)&0\dots 0\cr}.
$$
Let us note $\widetilde\PP'=\widetilde V-t{\widetilde{V}}^T$, where $\widetilde V=\pmatrix{N&M^T\cr M&\ell\cr}$. The computation above shows
that $\rk A(L';\Q)\ge\sum_{i=1}^n(d_i-1)$. The fact that the rank of $A(L';\Q)$ is equal to $\sum_{i=1}^n(d_i-1)$ can be proved by
(at least) two distinct methods. By a more subtle analysis of $V_\pm$, one can check that $\Delta_{L\m}=\det\widetilde\PP'\cdot\Delta'$
with some factor $\Delta'$; since $L\m$ is fibered, $\Delta_{L\m}\neq 0$ so $\det\widetilde\PP'\neq 0$ and
$\rk A(L';\Q)=\sum_{i=1}^n(d_i-1)$. A more conceptual proof goes as follows: $L'$ can be thought of as the result of the ``splicing'' of
$L\m$ with multilinks $L^{(1)}(\underline{m}^{(1)}),\dots,L^{(n)}(\underline{m}^{(n)})$ (see \cite{E-N,Neu}). It can be showed that
$\rk A(L^{(i)}(\underline{m}^{(i)}))=d_i-1$ for $i=1,\dots,n$, and that the rank of the Alexander module is additive under splicing
(see \cite[Theorem 4.3.1 and Proposition 3.2.4]{cim}). Since $L\m$ is fibered, $\rk A(L\m)=0$ and we get the result.

As a consequence, $\widetilde\PP'$ is a presentation matrix of the torsion submodule of $A(L';\Q)$. Now, note that
$$
(H_G^T\oplus I_{n-1}){\widetilde{V}}^T=\pmatrix{H_G^T&0\cr 0&I_{n-1}\cr}\pmatrix{N^T&M^T\cr M&\ell^T\cr}=
\pmatrix{H_G^TN^T&H_G^TM^T\cr M&\ell^T\cr}.
$$
By the equations $(\ast)$, this is exactly the matrix $\widetilde V$. Hence, the torsion submodule of $A(L';\Q)$ is presented by
$$
\widetilde\PP'=\widetilde V-t{\widetilde{V}}^T=(H^T_G\oplus I_{n-1}){\widetilde{V}}^T-t{\widetilde{V}}^T=((H_G^T\oplus I_{n-1})-tI)
{\widetilde{V}}^T.
$$
Since $\det\widetilde\PP'\neq 0$, ${\widetilde{V}}^T$ is unimodular. Therefore, $\widetilde\PP'$ is equivalent as a presentation matrix to
$(H_G^T\oplus I_{n-1})-tI=(H_G^T-tI)\oplus (1-t)I_{n-1}$. This concludes the proof.
\endproof

\section{Application to the Alexander module of links at infinity}

In this paragraph, we use Propositions \ref{prop:split} and \ref{prop:L'} to give a closed formula for the Alexander module over
$\C[t^{\pm 1}]$ of the regular link at infinity $\LL=\LL(f,\infty)$ of any polynomial map $f\colon\C^2\to\C$. Given such an $f$, there
exists a fibered multilink with multiplicities $m_i\neq 0$ and fiber $F$ such that $\LL=\partial F$. Furthermore, this multilink is an iterated
torus multilink: it can be constructed by iterated cabling and connected sum operations from the unknot. Since the Alexander module
over $\C[t^{\pm 1}]$ of iterated torus fibered multilinks is known,
the result for $\LL$ will follow directly from Propositions \ref{prop:split} and \ref{prop:L'}.

To state our result, we must assume that the reader is familiar with splice diagrams (see \cite{E-N}). Recall that a splice diagram
representing a multilink $L\m$ is a tree $\Gamma$ decorated as follows:
\begin{itemlist}
\item{Some of its leaves (valency one vertices) are drawn as arrowheads and represent components of $L$; they are endowed
with the multiplicity $m_i$ of the corresponding component $L_i$ of $L$.}
\item{Each edge has an integer weight at any end where it meets a node (vertex of valency greater than one), and these edge-weights
around a fixed node are pairwise coprime.}
\end{itemlist}  
Associated to each non-arrowhead vertex $v$ of $\Gamma$ is a so-called ``virtual component'': this is the additional link component
that would be represented by a single arrow at that vertex $v$ with edge-weight $1$. Splice diagrams are very convenient to compute
linking numbers: given two vertices $v$ and $w$ of $\Gamma$, the linking number of the corresponding components (virtual or ``real'')
is the product of all the edge-weights adjacent to but not on the shortest path in $\Gamma$ connecting $v$ and $w$. 

General splice diagrams as described here encode graph multilinks (that is: multilinks in homology sphere with graph manifold exterior).
A multilink in $S^3$ is a graph multilink if and only if it is an iterated torus multilink, so the multilink associated with a polynomial
map is encoded by such a splice diagram. Furthermore, Eisenbud and Neumann succeeded in computing the Alexander module over $\C[t^{\pm 1}]$
of any fibered graph multilink $L\m$ from its splice diagram $\Gamma$. If $L\m$ has ``uniform twists'' (this is the case of the multilink
associated with a polynomial map), the result goes as follows.

Let us denote by $\N$ the set of nodes of $\Gamma$, by $\E$ the set of edges connecting two nodes and by $\V$ the set of non-arrowhead vertices of $\Gamma$. By
cutting an edge $E\in\E$ in two, one gets two splice diagrams representing two multilinks; let us denote by $d_E$ the greatest common
divisor of the linking numbers of these two multilinks with the vitual component corresponding to the middle of the edge $E$. For every
$v\in\V$, let $\delta_v$ denote its valency and $\underline{m}(v)$ the linking number of $L\m$ with the virtual component corresponding to $v$.
Finally, for every node $v\in\N$, let $d_v$ be the greatest common divisor of the $d_E$'s of edges $E\in\E$ which meet
$v$, and of all the $m_i$'s of arrowheads adjacent to $v$.

\begin{thm}[{\rm Eisenbud-Neumann \cite[Theorem 14.1]{E-N}}]\label{thm:C-fibered}
Let $L\m$ be a fibered graph multilink with monodromy $h$ and uniform twists, given by a splice diagram $\Gamma$. 
The Alexander module $A(L\m;\C)$ is determined by the following properties:
\begin{itemlist}
\item{The Jordan normal form of $h_\ast$ consists of $\,1\times 1$ and $\,2\times 2$ Jordan blocks.}
\item{The characteristic polynomial of $h_\ast$ is equal to
$$
\Delta(t)=(t^d-1)\prod_{v\in\V}(t^{|\underline{m}(v)|}-1)^{\delta_v-2}.
$$}
\item{The eigenvalues corresponding to the $\,2\times2$ Jordan blocks are the roots of
$$
\Delta'(t)=(t^d-1)\frac{\prod_{E\in\E}(t^{d_E}-1)}{\prod_{v\in\N}(t^{d_v}-1)}.
$$}
\end{itemlist}
\end{thm}

\noindent Let us now state and prove our final result.

\begin{thm}\label{thm:infinity}
Let $f\colon\C^2\to\C$ be a polynomial map with regular link at infinity $\LL=\LL(f,\infty)$. If $L\m=L(m_1,\dots,m_n)$ denotes the
multilink associated with $\LL$, let $d$ be the greatest common divisor of $m_1,\dots,m_n$, and
$d_i=\gcd(m_i,\sum_{j\neq i}m_j\lk(L_i,L_j))$ for $i=1,\dots,n$. Also, let
$\Delta(t)$ be the characteristic polynomial of the monodromy of $L\m$, and $\Delta'(t)$ the polynomial corresponding to the
$\,2\times 2$ Jordan blocks (as in Theorem \ref{thm:C-fibered}).
Then, the Alexander module $A(\LL;\C)$ of $\LL$ over $\C[t^{\pm 1}]$ is given by the following properties:
\begin{itemlist}
\item{The rank of $A(\LL;\C)$ is equal to $\sum_{i=1}^n(d_i-1)$.}
\item{The Jordan normal form of $\,t\,$ restricted to the torsion submodule of $A(\LL;\C)$ consists of $\,1\times 1$ and $\,2\times 2$ Jordan blocks.}
\item{The order ideal of the torsion submodule of $A(\LL;\C)$ is generated by
$$
\widetilde{\Delta}(t)=(t-1)^{n-1}\frac{(t^d-1)\Delta(t)}{\prod_{i=1}^n(t^{d_i}-1)}\,.
$$}
\item{The eigenvalues corresponding to the $\,2\times2$ Jordan blocks are the roots of $\Delta'(t)$.}
\end{itemlist}
\end{thm}
\proof
The regular link at infinity $\LL$ is given by the boundary $\partial F$ of the fiber of $L\m$, which has non-zero multiplicities.
By Propositions \ref{prop:split} and \ref{prop:L'},
$$
A(\LL;\C)\;=\;(A_G\otimes\C[t^{\pm 1}])\oplus(\C[t^{\pm 1}]/(t-1))^{n-1}\oplus(\C[t^{\pm 1}])^{\sum_{i=1}^n(d_i-1)},
$$
where $A(L\m;\C)=(A_G\otimes\C[t^{\pm 1}])\oplus\bigoplus_{i=1}^{n-1}\C[t^{\pm 1}]\Big/\Bigl(\frac{(t^{D_i}-1)(t^{d_{i+1}}-1)}{(t^{D_{i+1}}-1)}\Bigr)$.
Therefore, the rank of $A(\LL;\C)$ is $\sum_{i=1}^n(d_i-1)$ and the order ideal of its torsion submodule is generated by
$$
(t-1)^{n-1}\frac{\Delta(t)}{\prod_{i=1}^{n-1}\frac{(t^{D_i}-1)(t^{d_{i+1}}-1)}{(t^{D_{i+1}}-1)}}=
(t-1)^{n-1}\frac{(t^d-1)\Delta(t)}{\prod_{i=1}^n(t^{d_i}-1)}\;,
$$
since $D_1=d_1$ and $D_n=\gcd(d_1,\dots,d_n)=\gcd(m_1,\dots,m_n)=d$. Furthermore, $A_G\otimes\C[t^{\pm 1}]$ contributes to Jordan blocks of dimension
at most two (by Theorem \ref{thm:C-fibered}), and
$\bigoplus_{i=1}^{n-1}\C[t^{\pm 1}]\Big/\Bigl(\frac{(t^{D_i}-1)(t^{d_{i+1}}-1)}{(t^{D_{i+1}}-1)}\Bigr)$ to Jordan blocks of dimension
one, since the polynomial $\frac{(t^{D_i}-1)(t^{d_{i+1}}-1)}{(t^{D_{i+1}}-1)}$ has only simple roots. \endproof

We refer to \cite[$\S$ 5.6]{cim} for a different proof of this result. Note that Propositions \ref{prop:split} and \ref{prop:L'} give
the Alexander module $A(\LL;\Q)$ from the module $A(L\m;\Q)$. The problem is that a closed formula for the Alexander module over
$\Q[t^{\pm 1}]$ of a fibered graph multilink remains unknown.

Let us conclude this note with an example.
  
\begin{ex}
Let $p,q,r$ be positive integers with $\gcd(p,r)=1$ and $p<(q+1)r$. Consider the polynomial map $f\colon\C^2\to\C$ given by
$$
f(x,y)=\left(x^qy+1\right)^r+x^p\,.
$$
As described in \cite[p. 451]{Neu}, the associated multilink
$L\m=L(1,qr)$ is given by the following splice diagram.

\begin{figure}[h]
  \begin{center}
  \epsfig{figure=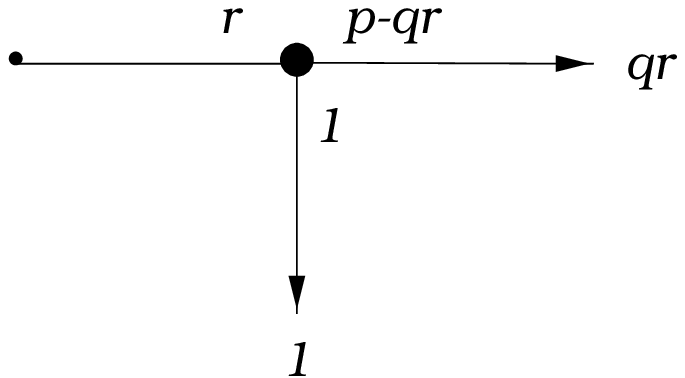,height=2.2cm}
  \end{center}
\end{figure}
\noindent Using Theorem \ref{thm:C-fibered}, one easily computes $\Delta(t)=(t-1)\frac{(t^{pr}-1)}{(t^p-1)}$ and $\Delta'(t)=1$. Hence
$$
A(L\m;\C)=\C[t^{\pm 1}]\Bigg/\left(\frac{(t-1)(t^{pr}-1)}{(t^p-1)}\right).
$$
Furthermore, $d_1=\gcd(1,qr^2)=1$ and $d_2=\gcd(qr,r)=r$. Therefore,
$$
A(\LL;\C)=\C[t^{\pm 1}]\Bigg/\left(\frac{(t-1)^2(t^{pr}-1)}{(t^p-1)(t^r-1)}\right)\oplus\left(\C[t^{\pm 1}]\right)^{r-1}.
$$
Note that $L\m$ is nothing but a torus multilink; on this simple example, it is possible to compute the Alexander modules
over $\Z[t^{\pm 1}]$. Using methods described in \cite{cim}, one can show that the Alexander module of $L\m$ is    
$\,\Z[t^{\pm 1}]\Big/\left(\frac{(t^{pr}-1)(t-1)}{(t^p-1)}\right)$, and that
$A(\LL)=\left(\Z[t^{\pm 1}]\right)^{r-1}\oplus\widetilde A(\LL)$, where $\widetilde A(\LL)$ is presented by the matrix
$$
\pmatrix{\frac{(t-1)(t^{pr}-1)}{(t^p-1)(t^r-1)}&q\cr0&q(t-1)\cr}.
$$
\end{ex}

\nonumsection{Acknowledgments}
I wish to thank Walter Neumann and Claude Weber for their encouragements, as well as Mathieu Baillif.

\nonumsection{References}

\end{document}